\documentclass[12pt]{article}

\usepackage{latexsym,amssymb}
\usepackage{graphicx}

\newtheorem{theorem}{Theorem}[section]

\newtheorem{example}[theorem]{Example}
\newtheorem{definition}[theorem]{Definition}

\newtheorem{proposition}[theorem]{Proposition}

\newtheorem{lemma}[theorem]{Lemma}
\newtheorem{corollary}[theorem]{Corollary}

\newtheorem{conjecture}[theorem]{Conjecture}

\newenvironment{proof}{\medskip\noindent{\it Proof.\ }}{\mbox{$\Box$}\medskip}

\begin{document}

\def\eqnsep{50pt}

\def\cM{\mathcal M}
\def\cN{\mathcal N}
\def\cO{\mathcal O}

\def\cplus{\ {{\ } \atop {+}}\ }
\def\cminus{\ {{\ } \atop {-}}\ }
\def\csum{\hbox{c}\hspace{-11pt}\sum}

\title{Restricted 3412-Avoiding Involutions:
Continued Fractions, Chebyshev Polynomials and Enumerations}
   
\author{
Eric S. Egge \\
Department of Mathematics \\
Gettysburg College\\
Gettysburg, PA  17325  USA \\
\\
eggee@member.ams.org \\
}

\maketitle

\begin{abstract}

Several authors have examined connections among restricted permutations, continued fractions, and Chebyshev polynomials of the second kind.
In this paper we prove analogues of these results for involutions which avoid 3412.
Our results include a recursive procedure for computing the generating function for involutions which avoid 3412 and any set of additional patterns.
We use our results to give enumerations and generating functions for involutions which avoid 3412 and various sets of additional patterns.
In many cases we express these generating functions in terms of Chebyshev polynomials of the second kind.

\medskip

\end{abstract}

\section{Introduction and Notation}

Let $S_n$ denote the set of permutations of $\{1, \ldots, n\}$, written in one-line notation, and suppose $\pi \in S_n$.
For all $i$, $1 \le i \le n$, we write $\pi(i)$ to denote the $i$th element of $\pi$.
We say $\pi$ is an {\it involution} whenever $\pi(\pi(i)) = i$ for all $i$, $1 \le i \le n$, and we write $I_n$ to denote the set of involutions in $S_n$.
Now suppose $\pi \in S_n$ and $\sigma \in S_k$.
We say a subsequence of $\pi$ has {\it type} $\sigma$ whenever it has all of the same pairwise comparisons as $\sigma$.
For example, the subsequence 2869 of the permutation 214538769 has type 1324.
We say $\pi$ {\it avoids} $\sigma$ whenever $\pi$ contains no subsequence of type $\sigma$.
For example, the permutation 214538769 avoids 312 and 2413, but it has 2586 as a subsequence so it does not avoid 1243.
If $\pi$ avoids $\sigma$ then $\sigma$ is sometimes called a {\it pattern} or a {\it forbidden subsequence} and $\pi$ is sometimes called a {\it restricted permutation} or a {\it pattern-avoiding permutation}.
In this paper we will be interested in permutations (and often only involutions) which avoid several patterns, so for any set $R$ of permutations we write $S_n(R)$ (resp. $I_n(R)$) to denote the set of permutations (resp. involutions) in $S_n$ which avoid every pattern in $R$ and we write $S(R)$ (resp. $I(R)$) to denote the set of all permutations (resp. involutions), including the empty permutation, which avoid every pattern in $R$.
When $R = \{\pi_1, \pi_2, \ldots, \pi_r\}$ we often write $S_n(R) = S_n(\pi_1, \pi_2, \dots, \pi_r)$, $I_n(R) = I_n(\pi_1, \pi_2, \dots, \pi_r)$, $S(R) = S(\pi_1, \pi_2, \ldots, \pi_r)$, and $I(R) = I(\pi_1, \pi_2, \dots, \pi_r)$.

As several authors have shown, generating functions for various subsets of $S(132)$ and $S(1243,2143)$ have close connections with continued fractions and Chebyshev polynomials of the second kind, and can often be computed recursively.
For example, Egge and Mansour \cite[Thm. 5.4]{EggeMansourSchroder} have shown that
\begin{equation}
\label{eqn:introcf12432143}
\sum_{\pi \in S(1243,2143)} \prod_{k \ge 1} x_k^{\tau_k(\pi)} = 1 + \frac{x_1}{\displaystyle
1 - x_1 - \frac{x_1 x_2}{\displaystyle
1 - x_1 x_2 - \frac{x_1 x_2^2 x_3}{\displaystyle
1 - x_1 x_2^2 x_3 - \frac{x_1 x_2^3 x_3^3 x_4}{\displaystyle
1 - x_1 x_2^3 x_3^3 x_4 - \cdots}}}},
\end{equation}
where $\tau_k(\pi)$ is the number of subsequences of type $12\ldots k$ in $\pi$.
Egge and Mansour have also shown \cite[Thm. 6.5]{EggeMansourSchroder} that
\begin{equation}
\label{eqn:introCheby}
\sum_{n=0}^\infty |S_n(1243, 2143, 12\ldots k)| x^n = 1 + \frac{\sqrt{x} U_{k-2}\left( \frac{1-x}{2 \sqrt{x}}\right)}{U_{k-1} \left( \frac{1-x}{2 \sqrt{x}} \right)},
\end{equation}
where $U_n(x)$ is the $n$th Chebyshev polynomial of the second kind, which may be defined by ${\displaystyle U_n(\cos t) = \frac{\sin((n+1) t)}{\sin t}}$.
Meanwhile, Krattenthaler has shown \cite[Theorem 3]{Krattenthaler} that
\begin{equation}
\label{eqn:intro12kr}
\sum_{\pi} x^{|\pi|} = \sum \prod_{i = 2}^b {{l_{i-1} + l_i - 1}\choose {l_i}} \frac{\left( U_{k-1}\left(\frac{1}{2 \sqrt{x}}\right)\right)^{l_1 - 1}}{\left(U_k\left( \frac{1}{2\sqrt{x}}\right)\right)^{l_1 + 1}} x^{\frac{1}{2}( l_1 - 1 ) + \sum\limits_{j=2}^b l_j}.
\end{equation}
Here the sum on the left is over all permutations in $S(132)$ which contain exactly $r$ subsequences of type $12\ldots k$, the quantity $|\pi|$ is the length of $\pi$, and the sum on the right is over all sequences $l_1, l_2, \ldots, l_b$ of nonnegative integers such that $\sum\limits_{i=1}^b l_i {{k+i-2}\choose{k-1}} = r$.
Finally, Mansour and Vainshtein have given \cite[Thm. 2.1]{MansourVainshtein2} the following recursive formula for $f_\pi(x) = \sum\limits_{n=0}^\infty |S_n(132,\pi)| x^n$, which makes it possible to compute $f_\pi(x)$ for any $\pi \in S(132)$.
\begin{equation}
\label{eqn:introFpirecurrence}
f_\pi(x) = 1 + x \sum_{j=0}^r \left( f_{\pi^j}(x) - f_{\pi^{j-1}}(x) \right) f_{\sigma^j}(x)
\end{equation}
Here $\pi^{j-1}$, $\pi^j$, and $\sigma^j$ are the types of certain subsequences of $\pi$.
For other results involving $S(132)$ and continued fractions or Chebyshev polynomials, see \cite{MansourVainshtein5} and the references therein.
For other results involving $S(1243,2143)$ and continued fractions or Chebyshev polynomials, see \cite{EggeMansourSchroder} and \cite{ARe}.

Involutions which avoid 3412 are known to have many properties which are analogous to properties of permutations which avoid 132.
For instance, it is well known that $|S_n(132)| = C_n$ for all $n \ge 0$, where $C_n$ is the $n$th Catalan number, which may be defined by $C_0 = 1$ and 
$$C_n = \sum_{i=1}^n C_{i-1} C_{n-i} \hspace{30pt} (n \ge 1).$$
(The Catalan number $C_n$ may also be defined by $C_n = \frac{1}{n+1} {{2n} \choose {n}}$.)
As a result, for all $n \ge 0$, the set $S_n(132)$ is in bijection with the set of {\it Dyck paths}.
These are the lattice paths from $(0,0)$ to $(2n,0)$ which contain only up $(1,1)$ and down $(1,-1)$ steps and which do not pass below the line $y = 0$.
Guibert \cite[Rem. 4.28]{G} has shown that $|I_n(3412)| = M_n$ for all $n \ge 0$, where $M_n$ is the $n$th Motzkin number, which may be defined by $M_0 = 1$ and
$$M_n = M_{n-1} + \sum_{i=2}^{n} M_{i-2} M_{n-i} \hspace{30pt} (n \ge 1).$$
As a result, for all $n \ge 0$, the set $I_n(3412)$ is in bijection with the set $\cM_n$ of {\it Motzkin paths}.
These are the lattice paths from $(0,0)$ to $(n,0)$ which contain only up $(1,1)$, down $(1,-1)$, and level $(1,0)$ steps and which do not pass below the line $y = 0$.
We write $\cM$ to denote the set of all Motzkin paths, including the empty path.
(For a partial list of other combinatorial objects counted by the Motzkin numbers, see \cite[pp. 238--9]{StanleyVol2}.)

Motivated by the parallels among $S(132)$, $S(1243, 2143)$, and $I(3412)$, in this paper we prove analogues of (\ref{eqn:introcf12432143}), (\ref{eqn:introCheby}), (\ref{eqn:intro12kr}), (\ref{eqn:introFpirecurrence}), and several similar results for $I(3412)$.
We begin with some results concerning $I(3412)$ and continued fractions.
We first define statistics $\tau_k$, $k \ge 1$, on $\cM$ and $I(3412)$.
On $I(3412)$, the statistic $\tau_k$ is simply the number of subsequences of type $k\ldots 21$.
On $\cM$, the statistic $\tau_k$ is a sum of binomial coefficients over the steps in the path.
We then give a simple bijection $\varphi : \cM \rightarrow I(3412)$ with the property that $\tau_k(\varphi(\pi)) = \tau_k(\pi)$ for all $k \ge 1$ and all $\pi \in \cM$.
Using $\varphi$ and a result of Flajolet, we prove the following analogue of (\ref{eqn:introcf12432143}).
\begin{equation}
\label{eqn:introIcf}
\sum_{\pi \in I(3412)} \prod_{k \ge 1} x_k^{\tau_k(\pi)} = 
\frac{1}{\displaystyle 1 - x_1 - \frac{x_1^2 x_2}{\displaystyle
1 - x_1 x_2^2 x_3 - \frac{x_1^2 x_2^5 x_3^4 x_4}{\displaystyle
1 - x_1 x_2^4 x_3^6 x_4^4 x_5 - \cdots}}}
\end{equation}
Here the $n$th numerator is $\prod\limits_{i=1}^{2n} x_i^{{{2n-2}\choose{i-1}} + {{2n-1}\choose{i-1}}}$ and the $n$th denominator is $1 - \prod\limits_{i=1}^{2n+1} x_i^{{{2n}\choose{i-1}}}$.
By specializing the $x_i$s in (\ref{eqn:introIcf}) appropriately, we obtain continued fraction expansions of the generating functions for several other statistics on $I_n(3412)$, including inversions, left-to-right maxima, and right-to-left minima.

We then turn our attention to analogues of (\ref{eqn:introCheby}), (\ref{eqn:intro12kr}), and (\ref{eqn:introFpirecurrence}).
We first use $\varphi$ and some well-known results concerning lattice paths to find the generating function for the involutions in $I(3412)$ which contain exactly $r$ subsequences of type $k\ldots 21$.
This generating function is an analogue of the generating function in (\ref{eqn:intro12kr}), and we express it in terms of Chebyshev polynomials.
We then find a recurrence relation for $F_T(x) = \sum\limits_{n=0}^\infty |I_n(3412,T)| x^n$ which enables us to compute $F_T(x)$ for any set $T$ of permutations. 
As a special case we have
\begin{equation}
\label{eqn:introIFpirecurrence}
F_\pi(x) = 1 + x F_\beta(x) + x^2 \sum_{i=1}^k \left( F_{\overline{\alpha_1 | \ldots | \alpha_i}}(x) - F_{\overline{\alpha_1 | \ldots | \alpha_{i-1}}}(x)\right) F_{\alpha_i | \ldots | \alpha_k}(x)
\end{equation}
for any permutation $\pi$.
Here the various subscripts of $F$ on the right are the types of certain subsequences of $\pi$.
This result is an analogue of (\ref{eqn:introFpirecurrence}).

Next we use (\ref{eqn:introIFpirecurrence}) to compute $F_\pi(x)$ for various $\pi$.
For example, we show that for all $k \ge 1$,
\begin{equation}
F_{2k \ldots 21}(x) = \frac{U_{k-1}\left(\frac{1-x}{2x}\right)}{x U_k\left(\frac{1-x}{2x}\right)}
\end{equation}
and
\begin{equation}
F_{2k-1 \ldots 21}(x) = \frac{U_{k-1}\left(\frac{1-x}{2x}\right) + U_{k-2}\left(\frac{1-x}{2x}\right)}{x \left( U_{k}\left(\frac{1-x}{2x}\right) + U_{k-1}\left(\frac{1-x}{2x}\right)\right)}.
\end{equation}
These results are analogues of (\ref{eqn:introCheby}). 
We also show that
\begin{equation}
F_{[k,l]}(x) = F_{k+l\ k+l-1\ldots 21}(x),
\end{equation}
where $[l_1,l_2,\ldots,l_m]$ is the layered permutation given by
$$[l_1, l_2, \ldots, l_m] = l_1, l_1-1, \ldots, 1, l_2 + l_1, l_2+l_1-1, \ldots, l_1+1, \ldots, \sum_{i=1}^m l_i, \sum_{i=1}^m l_i-1, \ldots, n-l_m+1.$$
For $m \ge 2$ the generating function $F_{[l_1,\ldots,l_m]}(x)$ does not reduce quite as nicely as it does when $m = 2$.
Nevertheless, we conjecture that $F_{[l_1,\ldots,l_m]}(x)$ is symmetric in $l_1,\ldots, l_m$ for all $m \ge 1$ and all $l_1,\ldots,l_m \ge 1$.
This conjecture has been verified for $m \le 4$ and $l_i \le 20$, as well as for $m = 5$ and $l_i \le 11$, using a Maple program.

In the next two sections we give enumerations and generating functions for various sets of involutions in $I(3412)$.
For instance, we show that for all $n \ge 0$ and all $k \ge 1$,
$$|I_n(3412, 4231, k+2\ k+1\ \ldots 4321)| = F_{k+1,n+1}.$$
Here $F_{k,n}$ is the $k$-generalized Fibonacci number, defined by $F_{k,n} = 0$ for $n \le 0$, $F_{k,1} = 1$, and $F_{k,n} = \sum\limits_{i=1}^k F_{k,n-i}$ for all $n \ge 1$.
We also give enumerations or generating functions for $I_n(3412, \sigma)$ for all $\sigma$ of length five or less.

We conclude the paper by mentioning several directions for future research.

\section{Two Families of Statistics and a Bijection}

In this section we give two infinite families of statistics, one on $\cM$ and the other on $I(3412)$.
We then give a simple bijection between $\cM_n$ and $I_n(3412)$ which relates our two families of statistics.
We begin by recalling the recursive structure of $\cM$.

\begin{definition}
For any Motzkin paths $\pi_1$ and $\pi_2$, we write
$$\pi_1 * \pi_2 = U \pi_1 D \pi_2.$$
\end{definition}

\begin{proposition}
\label{prop:pathmap}
\renewcommand\labelenumi{{\upshape (\roman{enumi}) }}
\begin{enumerate}
\item
For all $n \ge 1$, the map
$$
\begin{array}{ccc}
\cM_{n-1} &\longrightarrow& \cM_n \\
\pi &\mapsto& L \pi \\
\end{array}
$$
is a bijection between $\cM_{n-1}$ and the set of Motzkin paths in $\cM_n$ which begin with a level step.
\item
For all $n \ge 2$ and all $j$ such that $2 \le j \le n$, the map
$$
\begin{array}{ccc}
\cM_{j-2} \times \cM_{n-j} &\longrightarrow& \cM_n \\
(\pi_1,\pi_2) &\mapsto& \pi_1*\pi_2 \\
\end{array}
$$
is a bijection between $\cM_{j-2} \times \cM_{n-j}$ and the set of Motzkin paths in $\cM_n$ which begin with an up step and first return to the $x$-axis at $x=j$.
\end{enumerate}
\end{proposition}

Next we define our statistics on Motzkin paths.

\begin{definition}
Suppose $\pi$ is a Motzkin path.
For any step $s \in \pi$, we write $ht(s)$ to denote the {\em height} of $s$, which is the $x$-coordinate of the left-most point of $s$.
For any $k$ we write 
$$\tau_k(\pi) = \sum_{s \in \pi} {{2 ht(s)} \choose {k-1}} + \sum_{s \in \pi} {{2 ht(s) - 1} \choose {k-1}},$$
where the first sum on the right is over all up and level steps in $\pi$ and the second sum on the right is over all down steps in $\pi$.
Here we use the convention that ${{n} \choose {k}} = 0$ whenever $k < 0$ or $k > n$.
\end{definition}

\begin{example}
If $\pi = ULUUDLDDLUD$ then $\tau_1(\pi) = 11$, $\tau_2(\pi) = 22$, $\tau_3(\pi) = 27$, $\tau_4(\pi) = 19$, $\tau_5(\pi) = 7$, $\tau_6(\pi) = 1$, and $\tau_k(\pi) = 0$ for all other $k$.
\end{example}

As we show next, the statistics $\tau_k$ are compatible with the recursive structure of $\cM$.

\begin{proposition}
\label{prop:pathtauk}
\renewcommand\labelenumi{{\upshape (\roman{enumi}) }}
\begin{enumerate}
\item
For all $n \ge 1$, all $k$, and all $\pi \in \cM_{n-1}$,
\begin{equation}
\tau_k(L \pi) = {{0} \choose {k-1}} + \tau_k(\pi).
\end{equation}
\item
For all $n \ge 2$, all $j$ such that $2 \le j \le n$, all $k$, all $\pi_1 \in \cM_{j-2}$, and all $\pi_2 \in \cM_{n-j}$,
\begin{equation}
\tau_k(\pi_1 * \pi_2) = {{0} \choose {k-1}} + {{1} \choose {k-1}} + \tau_{k-2}(\pi_1) + 2 \tau_{k-1}(\pi_1) + \tau_k(\pi_1) + \tau_k(\pi_2).
\end{equation}
\end{enumerate} 
\end{proposition}
\begin{proof}
(i)
This is immediate from the definition of $\tau_k$.

(ii)
By definition of $\tau_k$ and $\pi_1 * \pi_2$ we have
$$\tau_k(\pi_1 * \pi_2) = {{0} \choose {k-1}} + \sum_{s \in \pi_1} {{2 ht(s) + 2} \choose {k-1}} + \sum_{s \in \pi_1} {{2 ht(s) + 1} \choose {k-1}} + {{1} \choose {k-1}} + \tau_k(\pi_2).$$
Here the first sum on the right is over all up and level steps in $\pi_1$ and the second sum on the right is over all down steps in $\pi_1$.
Now (ii) follows from the fact that for all $n \ge 0$ and all $k$, 
$${{n+2} \choose {k-1}} = {{n} \choose {k-3}} + 2 {{n} \choose {k-2}} + {{n} \choose {k-1}}.$$
\end{proof}

We now turn our attention to $I(3412)$, beginning with its recursive structure.
We start with some notation.

\begin{definition}
For any permutation $\pi$, we write $|\pi|$ to denote the length of $\pi$.
For any integer $n$ and any sequence $\sigma$ of integers we write $\sigma^{+n}$ to denote the sequence obtained by adding $n$ to every entry of $\sigma$.
\end{definition}

\begin{definition}
For any $\pi_1, \pi_2 \in I(3412)$ we write
$$\pi_1 * \pi_2 = |\pi_1|+2\ \pi_1^{+1}\ 1\ \pi_2^{+|\pi_1|+2}.$$
\end{definition}

We have now defined $*$ on both $\cM \times \cM$ and $I(3412) \times I(3412)$, but it will always be clear from the context which definition is intended.

\begin{proposition}
\label{prop:permmap}
\renewcommand\labelenumi{{\upshape (\roman{enumi}) }}
\begin{enumerate}
\item
For all $n \ge 1$, the map
$$
\begin{array}{ccc}
I_{n-1}(3412) &\longrightarrow& I_n(3412) \\
\pi & \mapsto & 1\ \pi^{+1}
\end{array}
$$
is a bijection between $I_{n-1}(3412)$ and the set of involutions in $I_n(3412)$ which begin with 1.
\item
For all $n \ge 0$ and all $j$ such that $2 \le j \le n$, the map
$$
\begin{array}{ccc}
I_{j-2}(3412) \times I_{n-j}(3412) &\longrightarrow& I_n(3412) \\
(\pi_1,\pi_2) &\mapsto& \pi_1 * \pi_2 \\
\end{array}
$$
is a bijection between $I_{j-2}(3412) \times I_{n-j}(3412)$ and the set of involutions in $I_n(3412)$ which begin with $j$.
\end{enumerate}
\end{proposition}
\begin{proof}
(i)
This is immediate.

(ii)
It is sufficient to show that the image of the given map is contained in $I_n(3412)$ and that the map is one-to-one and onto.

To see that the image of the given map is contained in $I_n(3412)$, suppose $\pi_1 \in I_{j-2}(3412)$ and $\pi_2 \in I_{n-j}(3412)$.
It is routine to verify that $\pi_1 * \pi_2$ is an involution, so we show it avoids 3412.
Suppose by way of contradiction that $abcd$ is a subsequence of $\pi_1 * \pi_2$ of type 3412.
Since $a > d$, $abcd$ must be entirely contained in $j, \pi_1^{+1}, 1$ or $\pi_2^{+|\pi_1|+2}$.
But $\pi_2$ avoids 3412, so $abcd$ is contained in $j, \pi_1^{+1}, 1$.
Since $a<b$, we must have $a \neq j$, and since $c<d$, we must also have $d \neq 1$, so $abcd$ is contained in $\pi_1^{+1}$.
But this contradicts our assumption that $\pi_1$ avoids $3412$, so $\pi_1 * \pi_2$ avoids 3412.

To see that the given map is one-to-one and onto, suppose $\pi \in I_n(3412)$ begins with $j$.
Since $\pi$ is an involution, we must also have $\pi(j) = 1$.
Moreover, if there exists $k > j$ such that $k$ appears between $j$ and 1 in $\pi$ then $\pi(k) \neq k$ and the subsequence $j, k, 1, \pi(k)$ has type 3412.
Since $\pi$ avoids 3412, the entries between $j$ and 1 in $\pi$ are $2, \ldots, j-1$.
It follows that there exist unique involutions $\pi_1 \in I_{j-2}(3412)$ and $\pi_2 \in I_{n-j}(3412)$ such that $\pi = \pi_1 * \pi_2$, so the given map is one-to-one and onto, as desired.
\end{proof}

We now define our statistics on $I(3412)$.

\begin{definition}
For any positive integer $k$ and any permutation $\pi$, we write $\tau_k(\pi)$ to denote the number of decreasing subsequences of length $k$ in $\pi$.
For notational convenience we set $\tau_k(\pi) = 0$ for all $k \le 0$.
\end{definition}

We have now defined $\tau_k$ on both $\cM$ and $I(3412)$, but it will always be clear from the context which definition is intended.
As we show next, the statistics $\tau_k$ are compatible with the recursive structure of $I(3412)$.

\begin{proposition}
\label{prop:permtauk}
\renewcommand\labelenumi{{\upshape (\roman{enumi}) }}
\begin{enumerate}
\item
For all $n \ge 1$, all $k$, and all $\pi \in I_{n-1}(3412)$ we have 
\begin{equation}
\label{eqn:taukDperm}
\tau_k(1\ \pi^{+1}) = {{0} \choose {k-1}} + \tau_k(\pi).
\end{equation}
\item
For all $n \ge 0$, all $j$ such that $2 \le j \le n$, all $k$, all $\pi_1 \in I_{j-2}(3412)$, and all $\pi_2 \in I_{n-j}(3412)$ we have
\begin{equation}
\label{eqn:tauk12perm}
\tau_k(\pi_1 * \pi_2) = {{0} \choose {k-1}} + {{1} \choose {k-1}} + \tau_{k-2}(\pi_1) + 2 \tau_{k-1}(\pi_1) + \tau_k(\pi_1) + \tau_k(\pi_2).
\end{equation}
\end{enumerate}
\end{proposition}
\begin{proof}
(i)
This is immediate from the definition of $\tau_k$.

(ii)
Observe that every decreasing subsequence of length at least two in $\pi_1 * \pi_2$ is entirely contained in either $|\pi_1| + 2\ \pi_1^{+1}\ 1$ or $\pi_2^{+|\pi_1|+2}$.
With this in mind, (ii) is immediate from the definition of $\pi_1 * \pi_2$.
\end{proof}

Next we introduce a bijection between $\cM_n$ and $I_n(3412)$ which is compatible with the statistics $\tau_k$ and the recursive structures of $\cM$ and $I(3412)$.

\begin{definition}
For any $\pi \in \cM_n$, we write $\varphi(\pi)$ to denote the permutation obtained as follows.
Number the steps in $\pi$ from left to right with $1,2,\ldots,n$.
For each up step at height $k$, find the first down step at height $k+1$ to its right and switch the labels of the two steps.
Then $\varphi(\pi)$ is the involution obtained by reading the resulting labels from left to right.
\end{definition}

\begin{example}
If $\pi = U L U D D L U D$ then $\varphi(\pi) = 52431687$.
\end{example}

\noindent
{\bf Remark}
The map $\varphi$ also appears in \cite[Rem. 4.28]{G}.

\begin{proposition}
\label{prop:varphibijection}
For all $n \ge 0$, the map $\varphi$ is a bijection between $I_n(3412)$ and $\cM_n$ such that $\tau_k(\varphi(\pi)) = \tau_k(\pi)$ for all $k$ and all $\pi \in \cM$.
\end{proposition}
\begin{proof}
Observe that if $\pi_1$ and $\pi_2$ are Motzkin paths then $\varphi(L \pi_1) = 1 \varphi(\pi)^{+1}$ and $\varphi(\pi_1 * \pi_2) = \varphi(\pi_1) * \varphi(\pi_2)$.
Arguing by induction on $n$, the result now follows from Propositions \ref{prop:pathmap}, \ref{prop:pathtauk}, \ref{prop:permmap}, and \ref{prop:permtauk}.
\end{proof}

For all $\pi \in \cM_n$, let $\pi^{rc}$ denote the Motzkin path in $\cM_n$ obtained from $\pi$ by reversing the order of the steps and switching the up and down steps.
For example, if $\pi = ULUUDLUUDDLLDDLUD$ then $\pi^{rc} = UDLUULLUUDDLUDDLD$.
Geometrically, $\pi^{rc}$ is the path obtained by reflecting $\pi$ over the line $x = \frac{n}{2}$.
We conclude this section by showing that the map on $I_n(3412)$ which corresponds via $\varphi$ to the map $\pi \mapsto \pi^{rc}$ is the well-known reverse complement map.

\begin{proposition}
\label{prop:reflectreversecomplement}
For any $\pi \in \cM$, let $\pi^{rc}$ be as in the paragraph above.
For any $\pi \in I_n(3412)$, let $\pi^{rc}$ denote the permutation obtained by reversing the order of the entries in $\pi$ and then replacing each entry $i$ with $n+1-i$.
Then $\varphi(\pi^{rc}) = \varphi(\pi)^{rc}$ for all $\pi \in \cM$.
\end{proposition}
\begin{proof}
For notational convenience we set $\alpha = \varphi(\pi)$, $\beta = \varphi(\pi^{rc})$, and $n = |\pi|$.

Observe that for any $i$, $1 \le i \le n$, we have $\alpha(i) = i$ if and only if the $i$th step of $\pi$ is a level step.
This happens if and only if the $n+1-i$th step of $\pi^{rc}$ is a level step, which happens if and only if $\beta(n+1-i) = n+1-i$.
It follows that $\alpha^{rc}$ and $\beta$ have the same fixed points.

Similarly, for all $i,j$ such that $1 \le i < j \le n$ we have $\alpha(i) = j$ and $\alpha(j) = i$ if and only if the $i$th step of $\pi$ is an up step and the $j$th step of $\pi$ is the corresponding down step.
This happens if and only if the $n+1-i$th step of $\pi^{rc}$ is a down step and the $n+1-j$th step of $\pi$ is the corresponding up step, which happens if and only if $\beta(n+1-i) = n+1-j$ and $\beta(n+1-j) = n+1-i$.
It follows that $\alpha^{rc}$ and $\beta$ have the same 2-cycles.

Combining these observations, we find that the entries of $\alpha^{rc}$ are the same as the entries of $\beta$, so $\alpha^{rc} = \beta$, as desired.
\end{proof}

\section{Continued Fractions}

In this section we combine our bijection $\varphi$ with a result of Flajolet to express the generating function for $I(3412)$ with respect to $\tau_k$, $k \ge 1$, as a continued fraction.
By specializing the indeterminates in this result appropriately, we also obtain continued fraction expansions of the generating functions for several other statistics on $I_n(3412)$, including inversions, left-to-right maxima, and right-to-left minima.
We begin by setting some notation.

\begin{definition}
For any given expressions $a_i$ $(i \ge 0)$ and $b_i$ $(i \ge 0)$ we write
$$\frac{a_0}{b_0} \cplus \frac{a_1}{b_1} \cplus \frac{a_2}{b_2} \cplus \frac{a_3}{b_3} \cplus \ldots$$
to denote the infinite continued fraction
$$
\frac{a_0}{\displaystyle
b_0 + \frac{a_1}{\displaystyle
b_1 + \frac{a_2}{\displaystyle
b_2 + \frac{a_3}{\displaystyle 
b_3 + \frac{a_4}{\displaystyle
b_4 + \cdots}}}}}.
$$
We use the corresponding notation for finite continued fractions.
\end{definition}

We now recall the relevant specialization of Flajolet's result.

\begin{theorem}
(Flajolet \cite[Theorem 1]{Flajolet})
For all $i \ge 1$, let $x_i$ denote an indeterminate.
Then we have
\begin{equation}
\label{eqn:pathcf}
\sum_{\pi \in \cM} \prod_{k \ge 1} x_k^{\tau_k(\pi)} = \frac{1}{1-x_1} \cminus \frac{x_1^2 x_2}{1 - x_1 x_2^2 x_3} \cminus \frac{x_1^2 x_2^5 x_3^4 x_4}{1 - x_1 x_2^4 x_3^6 x_4^4 x_5} \cminus \cdots \cminus \frac{\prod\limits_{i=1}^{2n} x_i^{{{2n-2} \choose {i-1}} + {{2n-1} \choose {i-1}}}}{1 - \prod\limits_{i=1}^{2n+1} x_i^{{{2n} \choose {i-1}}}} \cminus \cdots.
\end{equation}
\end{theorem}

Combining this with Proposition \ref{prop:varphibijection}, we obtain the following result.

\begin{theorem}
\label{thm:permcf}
For all $i \ge 1$, let $x_i$ denote an indeterminate.
Then we have
\begin{equation}
\label{eqn:permcf}
\sum_{\pi \in I(3412)} \prod_{k \ge 1} x_k^{\tau_k(\pi)} = \frac{1}{1-x_1} \cminus \frac{x_1^2 x_2}{1 - x_1 x_2^2 x_3} \cminus \frac{x_1^2 x_2^5 x_3^4 x_4}{1 - x_1 x_2^4 x_3^6 x_4^4 x_5} \cminus \cdots \cminus \frac{\prod\limits_{i=1}^{2n} x_i^{{{2n-2} \choose {i-1}} + {{2n-1} \choose {i-1}}}}{1 - \prod\limits_{i=1}^{2n+1} x_i^{{{2n} \choose {i-1}}}} \cminus \cdots.
\end{equation}
\end{theorem}

Using (\ref{eqn:permcf}) we can express the generating function for $I_n(3412, k\ldots 21)$ as a finite continued fraction.

\begin{corollary}
For all $k \ge 1$ we have
\begin{equation}
\label{eqn:avoid2k1}
\sum_{n=0}^\infty |I_n(3412, 2k \ldots 2 1)| x^n = \frac{1}{1-x} \cminus \underbrace{\frac{x^2}{1-x} \cminus \frac{x^2}{1-x} \cminus \cdots \cminus \frac{x^2}{1-x}}_{k-1\ terms}
\end{equation}
and
\begin{equation}
\label{eqn:avoid2k+11}
\sum_{n=0}^\infty |I_n(3412, 2k+1\ldots 2 1)| x^n = \frac{1}{1-x} \cminus \underbrace{\frac{x^2}{1-x} \cminus \frac{x^2}{1-x} \cminus \cdots \cminus \frac{x^2}{1-x}}_{k-1\ terms} \cminus \frac{x^2}{1}.
\end{equation}
\end{corollary}
\begin{proof}
To prove (\ref{eqn:avoid2k1}), set $x_1 = x$, $x_i = 1$ for all $i$, $2 \le i < 2k$, and $x_i = 0$ for all $i \ge 2k$ in (\ref{eqn:permcf}).

The proof of (\ref{eqn:avoid2k+11}) is similar to the proof of (\ref{eqn:avoid2k1}).
\end{proof}

We can also use (\ref{eqn:permcf}) to express the generating function for $I(3412)$ with respect to various statistics as a continued fraction.

\begin{corollary}
For any permutation $\pi$, let $inv(\pi)$ denote the number of inversions in $\pi$.
Then
$$\sum_{\pi \in I(3412)} q^{inv(\pi)} x^{|\pi|} = \frac{1}{1-x} \cminus \frac{x^2 q}{1-x q^2} \cminus \frac{x^2 q^5}{1 - xq^4} \cminus \cdots \cminus \frac{x q^{4n-3}}{1 - x q^{2n}} \cminus \cdots.$$
\end{corollary}
\begin{proof}
In (\ref{eqn:permcf}), set $x_1 = x$, $x_2 = q$, and $x_i = 1$ for all $i \ge 3$.
\end{proof}

\begin{corollary}
For any permutation $\pi$, let $m(\pi)$ denote the number of nonempty decreasing subsequences in $\pi$.
Then
$$\sum_{\pi \in I(3412)} q^{m(\pi)} x^{|\pi|} = \frac{1}{1-x} \cminus \frac{x^2 q^3}{1-x q^4} \cminus \frac{x^2 q^{12}}{1-x q^{16}} \cminus \cdots \cminus \frac{x^2 q^{3 \cdot 4^{n-1}}}{1-x q^{4^n}} \cminus \cdots.$$
\end{corollary}
\begin{proof}
In (\ref{eqn:permcf}), set $x_1 = xq$ and $x_i = q$ for all $i \ge 2$.
\end{proof}

For our next application of (\ref{eqn:permcf}), recall that a {\em left-to-right maximum} in a permutation $\pi$ is an entry of $\pi$ which is greater than all of the entries to its left.
Similarly, a {\em right-to-left minimum} in $\pi$ is an entry of $\pi$ which is less than all of the entries to its right.
As we show next, if $\pi \in I(3412)$ then the number of left-to-right maxima and the number of right-to-left minima in $\pi$ can be expressed in terms of the statistics $\tau_k$.
Combining this with (\ref{eqn:permcf}), we obtain a continued fraction expansion of the generating function for $I_n(3412)$ with respect to left-to-right maxima or right-to-left minima.

\begin{proposition}
\label{prop:maxmin}
For any permutation $\pi$, let $lrmax(\pi)$ denote the number of left-to-right maxima in $\pi$ and let $rlmin(\pi)$ denote the number of right-to-left minima in $\pi$.
Then for all $\pi \in I(3412)$,
$$lrmax(\pi) = rlmin(\pi) = \sum_{k=1}^\infty (-1)^{k-1} \tau_k(\pi).$$
\end{proposition}
\begin{proof}
Set $\tau = \sum\limits_{i=1}^\infty (-1)^{i+1} \tau_i$ and use Proposition \ref{prop:permtauk} to find that for all $\pi_1, \pi_2 \in I(3412)$ we have $\tau(1 \pi_1^{+1}) = 1 + \tau(\pi_1)$ and $\tau(\pi_1 * \pi_2) = 1 + \tau(\pi_2)$.
It is routine to verify that the same relations hold when $\tau$ is replaced with $lrmax$ or $rlmin$.
Using Proposition \ref{prop:permmap}, the result now follows by induction on $|\pi|$.
\end{proof}

\begin{corollary}
We have
$$\sum_{\pi \in I(3412)} q^{lrmax(\pi)} x^{|\pi|} = \sum_{\pi \in I(3412)} q^{rlmin(\pi)} x^{|\pi|} = \frac{1}{1-x} \cminus \frac{x^2 q}{1-x} \cminus \frac{x^2}{1-x} \cminus \cdots \cminus \frac{x^2}{1-x} \cminus \cdots.$$
\end{corollary}
\begin{proof}
In (\ref{eqn:permcf}), set $x_1 = xq$ and $x_i = q^{(-1)^{i-1}}$ for all $i \ge 2$ and use Proposition \ref{prop:maxmin} to simplify the result.
\end{proof}

For our final application of (\ref{eqn:permcf}), recall that $i$ is a fixed point for a permutation $\pi$ whenever $\pi(i) = i$.
As we show next, if $\pi \in I(3412)$ then the number of fixed points in $\pi$ can be expressed in terms of the statistics $\tau_k$.
Combining this with (\ref{eqn:permcf}), we obtain a continued fraction expansion of the generating function for $I_n(3412)$ with respect to the number of fixed points.

\begin{proposition}
\label{prop:fixedpoints}
For any permutation $\pi$, let $fix(\pi)$ denote the number of fixed points in $\pi$.
Then for all $\pi \in I(3412)$,
$$fix(\pi) = \sum_{k=1}^\infty (-2)^{k-1} \tau_k(\pi).$$
\end{proposition}
\begin{proof}
This is similar to the proof of Proposition \ref{prop:maxmin}.
\end{proof}

\begin{corollary}
We have
$$\sum_{\pi \in I(3412)} q^{fix(\pi)} x^{|\pi|} = \frac{1}{1-xq} \cminus \frac{x^2}{1-xq} \cminus \frac{x^2}{1-xq} \cminus \cdots \cminus \frac{x^2}{1-xq} \cminus \cdots.$$
\end{corollary}
\begin{proof}
In (\ref{eqn:permcf}), set $x_1 = xq$ and $x_i = q^{(-2)^{i-1}}$ for all $i \ge 2$ and use Proposition \ref{prop:fixedpoints} to simplify the result.
\end{proof}

We now turn our attention to the question of which statistics on $I(3412)$ have generating functions which can be expressed as continued fractions like the one in (\ref{eqn:permcf}).
We begin by specifying which continued fractions we wish to consider.

By a {\em Motzkin continued fraction} we mean a continued fraction of the form
$$\frac{1}{1-m_0} \cminus \frac{m_0 m_1}{1-m_2} \cminus \frac{m_2 m_3}{1-m_4} \cminus \cdots \cminus \frac{m_{2n-2} m_{2n-1}}{1 - m_{2n}} \cminus \cdots,$$
where $m_i$ is a monic monomial in a given set of variables for all $i \ge 0$.
Observe that if $f_1, f_2, f_3, \ldots$ are (possibly infinite) linear combinations of the $\tau_k$s with the property that each $\tau_k$ appears in only finitely many $f_i$, then by specializing the $x_i$s appropriately in (\ref{eqn:permcf}) we can express the generating function
$$\sum_{\pi \in I(3412)} x^{|\pi|} \prod_{k \ge 1} q_k^{f_k(\pi)}$$
as a Motzkin continued fraction.
For example, when only $f_1$ is present, we have the following corollary of Theorem \ref{thm:permcf}.

\begin{corollary}
Let $\lambda_1, \lambda_2, \ldots$ denote nonnegative integers and let $f$ denote the statistic
$$f = \sum_{k \ge 1} \lambda_k \tau_k$$
on $I(3412)$.
Then
\begin{eqnarray*}
\lefteqn{\sum_{\pi \in I(3412)} q^{f(\pi)} x^{|\pi|} = \frac{1}{1-x q^{f(1)}} \cminus \frac{x^2 q^{f(21)}}{1-x q^{f(321) - f(21)}} \cminus} & & \\
& & \frac{x^2 q^{f(4321) - f(21)}}{1 - x q^{f(54321) - f(4321)}} \cminus \cdots \cminus \frac{x^2 q^{f(2n\ldots 2 1) - f(2n-2\ldots 2 1)}}{1 - x q^{f(2n+1\dots 21) - f(2n\ldots 2 1)}} \cminus \cdots.\\
\end{eqnarray*}
\end{corollary}
\begin{proof}
In (\ref{eqn:permcf}), set $x_1 = x q^{\lambda_1}$ and $x_i = q^{\lambda_i}$ for all $i \ge 2$ and use the fact that 
$$f(n \ldots 2 1) - f(n-1 \ldots 2 1) = \sum_{i=0}^{n-1} {{n-1} \choose {i}} \lambda_i$$
for all $n \ge 2$ to simplify the result.
\end{proof}

With the same arguments used to prove \cite[Thm. 2]{BCS} and \cite[Thm. 5.12]{EggeMansourSchroder}, one can also prove the following result.

\begin{theorem}
The set of Motzkin continued fractions is exactly the set of generating functions for countable families of statistics on $I(3412)$ in which each statistic is a (possibly infinite) linear combination of the $\tau_k$s and each $\tau_k$ appears in only finitely many statistics.
\end{theorem} 

\section{Involutions Avoiding 3412 and Containing $k \ldots 21$}
\label{sec:rk21}

In this section we use our bijection $\varphi$ to find the generating function for the permutations in $I(3412)$ which contain exactly $r$ decreasing subsequences of length $k$.
We express this generating function in terms of Chebyshev polynomials of the second kind, so we begin by recalling these polynomials.

\begin{definition}
For all $n \ge -1$, we write $U_n(x)$ to denote the {\em $n$th Chebyshev polynomial of the second kind}, which is defined by $U_{-1}(x) = 0$ and ${\displaystyle U_n(\cos t) = \frac{\sin ((n+1)t)}{\sin t}}$ for $n \ge 0$.
These polynomials satisfy
\begin{equation}
\label{eqn:Chebyshevrecurrence}
U_n(x) = 2x U_{n-1}(x) - U_{n-2}(x) \hspace{30pt} (n \ge 1).
\end{equation}
\end{definition}

We will find it useful to reformulate the recurrence in (\ref{eqn:Chebyshevrecurrence}), replacing $x$ with ${\displaystyle \frac{1-x}{2x}}$ to obtain
\begin{equation}
\label{eqn:ourChebyrecurrence}
x U_n\left( \frac{1-x}{2x}\right) = (1-x) U_{n-1}\left(\frac{1-x}{2x}\right) - x U_{n-2}\left(\frac{1-x}{2x}\right).
\end{equation}

Our main results in this section are the following, which give the promised generating function.

\begin{theorem}
\label{thm:biggfeven}
Fix $r \ge 1$, $k \ge 1$, and $b \ge 0$ such that
$$r < min\left( {{2k+2b+2} \choose {2k-1}}, {{2k+2b}\choose{2k-1}} + {{2k+2b+1}\choose{2k-1}}\right).$$
Then
\begin{equation}
\label{eqn:biggfeven}
\sum_{\pi} x^{|\pi|} = \sum \prod_{i=0}^b {{d_i + d_{i+1} + l_i - 1}\choose{d_{i+1}+l_i}} {{d_{i+1}+l_i}\choose{l_i}} \frac{\left(U_{k-1}\left(\frac{1-x}{2x}\right)\right)^{d_0-1}}{\left(U_k\left(\frac{1-x}{2x}\right)\right)^{d_0+1}} x^{-1-d_0+\sum\limits_{j=0}^b (2 d_i + l_i)}.
\end{equation}
Here the sum on the left is over all involutions in $I(3412)$ which contain exactly $r$ subsequences of type $2k\ldots 21$.
The sum on the right is over all sequences $d_0, \ldots, d_b$ and $l_0, \ldots, l_b$ of nonnegative integers such that
\begin{equation}
\label{eqn:bigreven}
r = \sum_{i=0}^b d_i \left( {{2k+2i-2}\choose{2k-1}} + {{2k+2i-1}\choose{2k-1}}\right) + \sum_{i=0}^b l_i {{2k+2i}\choose{2k-1}}.
\end{equation}
Throughout we adopt the convention that ${{a} \choose {0}} = 1$ and ${{a} \choose {-1}} = 0$ for any integer $a$.
\end{theorem}

\begin{theorem}
\label{thm:biggfodd}
Fix $r \ge 1$, $k \ge 1$, and $b \ge 0$ such that ${{2k+2b}\choose{2k}} \le r < {{2k+2b+2}\choose{2k}}$.
Then
\begin{equation}
\label{eqn:biggfodd}
\sum_{\pi} x^{|\pi|} = \sum \prod_{i=0}^b {{d_i+d_{i-1}+l_i-1}\choose{d_i+l_i}}{{d_i+l_i}\choose{l_i}} \frac{\left(U_k\left(\frac{1-x}{2x}\right)\right)^{d_0+l_0-1}}{\left(U_{k+1}\left(\frac{1-x}{2x}\right) + U_k\left(\frac{1-x}{2x}\right)\right)^{d_0+l_0+1}} x^{-1-d_0-l_0+\sum\limits_{j=0}^b (2d_j+l_j)}.
\end{equation}
Here the sum on the left is over all involutions in $I(3412)$ which contain exactly $r$ subsequences of type $2k+1\ldots 21$.
The sum on the right is over all sequences $d_0, \ldots, d_b$ and $l_0, \ldots, l_b$ of nonnegative integers such that
\begin{equation}
\label{eqn:bigrodd}
r = \sum_{i=0}^b d_i\left( {{2k+2i+1}\choose{2k}} + {{2k+2i}\choose{2k}}\right) + \sum_{i=0}^b l_i{{2k+2i}\choose{2k}}.
\end{equation}
For notational convenience we set $d_{-1} = 1$.
Throughout we adopt the convention that ${{a} \choose {0}} = 1$ and ${{a} \choose {-1}} = 0$ for any integer $a$.
\end{theorem}

To prove these theorems, we first need to set some notation and recall some preliminary results.  
We begin with some matrices which will prove useful.

\begin{definition}
For all $k \ge 0$ we write $A_k$ to denote the $k+1$ by $k+1$ tridiagonal matrix given by
$$A_k = \left( \matrix{x & x & 0 & 0 & 0 & \cdots & 0 & 0 & 0 & 0 \cr
x & x & x & 0 & 0 & \cdots & 0 & 0 & 0 & 0 \cr
0 & x & x & x & 0 & \cdots & 0 & 0 & 0 & 0 \cr
\vdots & & \ddots & \ddots & \ddots & & & & \vdots & \vdots \cr
\vdots & & & \ddots & \ddots & \ddots & & & \vdots & \vdots \cr
\vdots & & & & \ddots & \ddots & \ddots & & \vdots & \vdots \cr
\vdots & & & & & \ddots & \ddots & \ddots & \vdots & \vdots \cr
0 & 0 & 0 & 0 & 0 & \cdots & x & x & x & 0 \cr
0 & 0 & 0 & 0 & 0 & \cdots & 0 & x & x & x \cr
0 & 0 & 0 & 0 & 0 & \cdots & 0 & 0 & x & x}
\right).$$
We write $B_k$ to denote the $k+1$ by $k+1$ tridiagonal matrix obtained by replacing the entry in the lower right corner of $A_k$ with 0.
We write $C_k$ to denote the $k+1$ by $k+1$ tridiagonal matrix obtained by replacing the entry in the upper left corner of $A_k$ with 0.
\end{definition}

The matrices $A_k$, $B_k$, and $C_k$ are closely related to generating functions for various sets of Motzkin paths.
To describe this relationship, we let $\cM(r,s,k)$ denote the set of lattice paths involving only up $(1,1)$, down $(1,-1)$, and level $(1,0)$ steps which begin at a point at height $r$, $0 \le r \le k$, end at a point at height $s$, $0 \le s \le k$, and do not cross the lines $y = k$ and $y = 0$.
Similarly, we let $\cN(r,s,k)$ denote the set of lattice paths in $\cM(r,s,k)$ which do not have any level steps at height $k$, and we let $\cO(r,s,k)$ denote the set of lattice paths in $\cM(r,s,k)$ which do not have any level steps at height 0.
Modifying the proof of \cite[Thm. A2]{Krattenthaler} slightly, we find that
\begin{equation}
\label{eqn:Agf}
\sum_{\pi \in \cM(r,s,k)} x^{|\pi|} = \frac{(-1)^{r+s} \det(I - A_k; s, r)}{\det(I - A_k)},
\end{equation}
\begin{equation}
\label{eqn:Bgf}
\sum_{\pi \in \cN(r,s,k)} x^{|\pi|} = \frac{(-1)^{r+s} \det(I - B_k; s, r)}{\det(I - B_k)},
\end{equation}
and
\begin{equation}
\label{eqn:Cgf}
\sum_{\pi \in \cO(r,s,k)} x^{|\pi|} = \frac{(-1)^{r+s} \det(I - C_k; s, r)}{\det(I - C_k)}.
\end{equation}
Here $|\pi|$ is the number of steps in $\pi$, $I$ is the identity matrix of the appropriate size, and $\det(M; s, r)$ is the minor of the matrix $M$ in which the $s$th row and $r$th column of $M$ have been deleted.
The determinants in (\ref{eqn:Agf}), (\ref{eqn:Bgf}), and (\ref{eqn:Cgf}) can often be expressed in terms of Chebyshev polynomials of the second kind.
For instance, arguing by induction on $k$ we find that for all $k \ge 0$,
\begin{equation}
\label{eqn:Auk}
x^{k+1} U_{k+1}\left(\frac{1-x}{2x}\right) = \det(I - A_k)
\end{equation}
and
\begin{equation}
\label{eqn:Buk}
x^{k+1} \left( U_{k+1}\left(\frac{1-x}{2x}\right) + U_k\left(\frac{1-x}{2x}\right) \right) = \det(I - B_k) = \det(I-C_k).
\end{equation}

We now prove Theorem \ref{thm:biggfeven}.

\noindent
{\em Proof of Theorem \ref{thm:biggfeven}}.
First observe that in view of Proposition \ref{prop:varphibijection}, the generating function on the left side of (\ref{eqn:biggfeven}) is the generating function for the set of Motzkin paths from $(0,0)$ to $(n,0)$ for which $\tau_{2k}(\pi) = r$.
To compute this generating function, observe that every Motzkin path $\pi$ with $\tau_{2k}(\pi) = r$ can be constructed by the following procedure in exactly one way.
\begin{enumerate}
\item
Choose $d_0, \ldots, d_b$ and $l_0, \ldots, l_b$ such that (\ref{eqn:bigreven}) holds.
Construct a sequence of down and level steps which contains exactly $d_i$ down steps at height $k+i$ and $l_i$ level steps at height $k+i$ for $0 \le i \le b$ and which satisfies all of the following.
\begin{enumerate}
\item
The step immediately preceeding a step at height $j$ is either a down step at height $j+1$ or less or a level step at height $j$ or less.
\item
All steps after the last down step at height $j$ are at height $j-1$ or less.
\item
The sequence ends with a down step at height $k$.
\end{enumerate}
\item
If the first step is at height $k+j$, insert $j+1$ up steps before the first step.
Similarly, after each step except the last, insert enough up steps to reach the height of the next level or down step.
\item
After each down step at height $k$ except the last, insert an (possibly empty) upside-down Motzkin path of height at most $k-1$.
\item
Before the first step insert a path from height 0 to height $k-1$ which does not exceed height $k-1$.
\item
After the last step, insert a path from height $k-1$ to height 0 which does not exceed height $k-1$.
\end{enumerate}
Since the choice at each step is independent of the choices at the other steps, and since every sequence of choices results in a path of the type desired, the desired generating function is the product of the generating functions for each step.

To compute the generating function for step 1, suppose we have fixed $d_0, \ldots, d_b$ and $l_0, \ldots, l_b$;  then each of the resulting partial paths will have generating function $x^{\sum\limits_{j=0}^b (d_i+l_i)}$.
To count these paths, we construct them from the top down.
That is, we first arrange the $l_b$ level steps at height $k+b$;  there is one way to do this.
We then place the $d_b$ down steps at height $k+b$ so that one of these steps occurs after all of the diagonal steps.
There are ${{d_b + l_b - 1} \choose {l_b}}$ ways to do this.
We then place the $l_{b-1}$ level steps at height $k+b-1$ so that none of these steps immediately follows a diagonal step at height $k+b$.
There are ${{d_b + l_{b-1}}\choose{l_{b-1}}}$ ways to do this.
Proceeding in this fashion, we find that the generating function for step 1 is equal to
\begin{equation}
\label{eqn:step1gf}
\sum \prod_{i=0}^b {{d_i + d_{i+1} + l_i - 1}\choose{d_{i+1} + l_i}} {{d_{i+1} + l_i}\choose{l_i}} x^{\sum\limits_{j=0}^b (d_i + l_i)},
\end{equation}
where the sum on the left is over all sequences $d_0, \ldots, d_b$ and $l_0, \ldots, l_b$ of nonnegative integers which satisfy (\ref{eqn:bigreven}).
In the path obtained after step 2 there is exactly one up step for every down step so the generating function for step 2 is equal to 
\begin{equation}
\label{eqn:step2gf}
x^{\sum\limits_{j=0}^b d_i}.
\end{equation}
Using (\ref{eqn:Agf}) and (\ref{eqn:Auk}), we find that the generating function for step 3 is equal to 
\begin{equation}
\label{eqn:step3gf}
\left( \frac{U_{k-1}\left(\frac{1-x}{2x}\right)}{x U_k\left(\frac{1-x}{2x}\right)}\right)^{d_0-1}
\end{equation}
and the generating functions for steps 4 and 5 are both equal to
\begin{equation}
\label{eqn:step45gf}
\frac{1}{x U_k\left(\frac{1-x}{2x}\right)}.
\end{equation}
Taking the product of the quantities in (\ref{eqn:step1gf}), (\ref{eqn:step2gf}), (\ref{eqn:step3gf}) and the square of the quantity in (\ref{eqn:step45gf}), we obtain (\ref{eqn:biggfeven}), as desired.
$\Box$

The proof of Theorem \ref{thm:biggfodd} is similar to the proof of Theorem \ref{thm:biggfeven}, using (\ref{eqn:Bgf}), (\ref{eqn:Cgf}), and (\ref{eqn:Buk}).

Theorems \ref{thm:biggfeven} and \ref{thm:biggfodd} have several interesting special cases;  we give four of them here.

\begin{corollary}
For all $k \ge 1$,
$$\sum_{\pi} x^{|\pi|} = \frac{1}{\left(U_k\left(\frac{1-x}{2x}\right)\right)^2},$$
where the sum on the left is over all involutions in $I(3412)$ which contain exactly one subsequence of type $2k \ldots 21$.
\end{corollary}

\begin{corollary}
\label{cor:3211}
For all $k \ge 1$,
$$\sum_{\pi} x^{|\pi|} = \frac{1}{x \left( U_{k+1}\left(\frac{1-x}{2x}\right) + U_k\left(\frac{1-x}{2x}\right)\right)^2},$$
where the sum on the left is over all involutions in $I(3412)$ which contain exactly one subsequence of type $2k+1\ldots 21$.
\end{corollary}

\begin{corollary}
For all $k \ge 1$,
$$\sum_{\pi} x^{|\pi|} = \frac{x U_{k-1}\left(\frac{1-x}{2x}\right)}{\left( U_k\left(\frac{1-x}{2x}\right)\right)^3},$$
where the sum on the left is over all involutions in $I(3412)$ which contain exactly two subsequences of type $2k \ldots 21$.
\end{corollary}

\begin{corollary}
\label{cor:3212}
For all $k \ge 1$,
$$\sum_{\pi} x^{|\pi|} = \frac{U_k\left(\frac{1-x}{2x}\right)}{x \left( U_{k+1}\left(\frac{1-x}{2x}\right) + U_k\left(\frac{1-x}{2x}\right)\right)^3},$$
where the sum on the left is over all involutions in $I(3412)$ which contain exactly two subsequences of type $2k+1\ldots 21$.
\end{corollary}

We close this section by observing that for all $r \ge 1$, the number of involutions in $I_n(3412)$ which contain exactly $r$ subsequences of type 321 can always be expressed in terms of Fibonacci numbers, since $U_2\left(\frac{1-x}{2x}\right) + U_1\left(\frac{1-x}{2x}\right) = \frac{1-x-x^2}{x^2}$.
For example, it follows from Corollaries \ref{cor:3211} and \ref{cor:3212} that the number of involutions in $I_n(3412)$ which contain exactly one subsequence of type 321 is ${\displaystyle \frac{2(n-1)F_n - n F_{n-1}}{5}}$ and the number of involutions in $I_n(3412)$ which contain exactly two subsequences of type 321 is ${\displaystyle \frac{5n^2 - 9n}{25} F_{n+1} + \frac{-15n^2+29n-6}{50} F_n}$.

\section{Restricted 3412-Avoiding Involutions}

We now turn our attention to generating functions for involutions in $I(3412)$ which avoid a set of additional patterns.
For any set $T$ of permutations we write
$$F_T(x) = \sum_{n = 0}^\infty |I_n(3412, T)| x^n,$$
and we observe that $F_\emptyset(x) = 0$, $F_1(x) = 1$, and ${\displaystyle F_{12}(x) = F_{21}(x) = \frac{1}{1-x}}$.
In this section we give a recurrence relation which allows one to compute $F_T(x)$ for any $T$.
We begin with a method of decomposing permutations and a map on permutations.

\begin{definition}
Fix $n \ge 1$.
We call a permutation $\pi \in S_n$ {\em complete} whenever no initial segment of $\pi$ of length $k$, $1 \le k < n$, consists of the numbers $1, 2, \ldots, k$.
\end{definition}

Observe that for every permutation $\pi$ there exists a unique sequence of complete permutations $\alpha_1, \ldots, \alpha_k$ such that $\pi = \alpha_1 \alpha_2^{+|\alpha_1|} \alpha_3^{+|\alpha_1|+|\alpha_2|} \ldots \alpha_k^{+\sum\limits_{i=1}^{k-1} |\alpha_i|}.$
In this situation we often abbreviate $\pi = \alpha_1 | \alpha_2 | \ldots | \alpha_k$.

\begin{definition}
\label{defn:overline}
For any permutation $\pi$, we define $\overline{\pi}$ as follows.
\begin{enumerate}
\item
$\overline{\emptyset} = \emptyset$ and $\overline{1} = \emptyset$.
\item
If $|\pi| \ge 2$ and there exists a permutation $\sigma$ such that $\pi = |\pi| \sigma^{+1} 1$ then $\overline{\pi} = \sigma$.
\item
If $|\pi| \ge 2$, there exists a permutation $\sigma$ such that $\pi = |\pi| \sigma$, and $\sigma$ does not end with 1 then $\overline{\pi} = \sigma$.
\item
If $|\pi| \ge 2$, there exists a permutation $\sigma$ such that $\pi = \sigma^{+1} 1$, and $\pi$ does not begin with $|\pi|$ then $\overline{\pi} = \sigma$.
\item
If $|\pi| \ge 2$, $\pi$ does not begin with $|\pi|$, and $\pi$ does not end with 1 then $\overline{\pi} = \pi$.
\end{enumerate}
\end{definition}

Observe that if $\pi$ and $\sigma$ are permutations then $|\pi|+2\ \pi^{+1} 1$ avoids $\sigma$ if and only if $\pi$ avoids $\overline{\sigma}$.

In order to give our recurrence relation for $F_T(x)$, we need to set some additional notation.

\begin{definition}
\label{defn:setnotation}
Let $T = \{\pi_1,\ldots,\pi_m\}$ denote a set of permutations and fix complete permutations $\alpha^i_j$, $1 \le i \le m$, $1 \le j \le k_i$, such that $\pi_i = \alpha^i_1 | \ldots | \alpha^i_{k_i}$.
For all $i_1,\ldots, i_m$ such that $0 \le i_j \le k_j$, let $T^{left}_{i_1,\ldots,i_m} = \{\overline{\alpha^1_1| \ldots |\alpha^1_{i_1}}, \ldots, \overline{\alpha^m_1|\ldots|\alpha^m_{i_m}}\}$ and let $T^{right}_{i_1,\ldots,i_m} = \{\alpha^1_{i_1}|\ldots|\alpha^1_{k_1},\ldots,\alpha^m_{i_m}|\ldots|\alpha^m_{k_m}\}$.
For any subset $Y \subseteq \{1,\ldots,m\}$, set
\begin{equation}
T_Y = \bigcup_{j \in Y} \{\overline{\alpha^j_1|\ldots|\alpha^j_{i_j-1}}\} \bigcup_{j \not\in Y, 1 \le j \le m} \{\overline{\alpha^j_1|\ldots|\alpha^j_{i_j}}\}.
\end{equation}
\end{definition}

We now describe how to find the generating function for permutations which contain patterns in one set while avoiding patterns in another set.

\begin{lemma}
\label{lem:Y}
With reference to Definition \ref{defn:setnotation}, fix $i_1,\ldots, i_m$ such that $1 \le i_j \le k_j$.
Then the generating function for those permutations which contain every pattern in $T^{left}_{i_1-1,\ldots,i_m-1}$ and avoid every pattern in $T^{left}_{i_1,\ldots,i_m}$ is
\begin{equation}
\sum_{Y \subseteq \{1,2,\ldots,m\}} (-1)^{|Y|} F_{T_Y}(x).
\end{equation}
\end{lemma}
\begin{proof}
This follows by a routine inclusion-exclusion argument, since $S(3412,T_{Y_1}) \subseteq S(3412,T_{Y_2})$ whenever $Y_2 \subseteq Y_1$.
\end{proof}

We are now ready to give our recurrence relation for $F_T(x)$.

\begin{theorem}
\label{thm:FTrecurrence}
With reference to Definition \ref{defn:setnotation},
\begin{equation}
\label{eqn:FTrecurrence}
F_T(x) = 1 + x F_{\beta(T)}(x) + x^2 \sum_{i_1,\ldots,i_m = 1}^{k_1,\ldots,k_m} \left( \sum_{Y \subseteq \{1,2,\ldots,m\}} (-1)^{|Y|} F_{T_Y}(x) \right) F_{T^{right}_{i_1,\ldots,i_m}}(x).
\end{equation}
Here $\beta(\pi_i) = \pi_i$ if $\alpha^i_1 \neq 1$, $\beta(\pi_i) = \alpha^i_2| \ldots | \alpha_{k_i}^i$ if $\alpha^i_1 = 1$, and $\beta(T)$ is the set of permutations obtained by applying $\beta$ to every element of $T$.
\end{theorem}
\begin{proof}
The set $I(3412,T)$ can be partitioned into three sets:  the set $A_1$ containing only the empty permutation, the set $A_2$ of those involutions which begin with 1, and the set $A_3$ of those involutions which do not begin with 1.

The generating function for $A_1$ is 1.

In view of Proposition \ref{prop:permmap}(i), the generating function for $A_2$ is $x F_{\beta(T)}(x)$, where $\beta(\pi_i) = \pi_i$ if $\alpha^i_1 \neq 1$, $\beta(\pi_i) = \alpha^i_2| \ldots | \alpha_{k_i}^i$ if $\alpha^i_1 = 1$, and $\beta(T)$ is the set of permutations obtained by applying $\beta$ to every element of $T$.

To obtain the generating function for $A_3$, we first observe that in view of Proposition \ref{prop:permmap}(ii), all permutations in $A_3$ have the form $\sigma_1 * \sigma_2$.
Since each $\alpha^i_j$ is complete, if $\sigma_1 * \sigma_2$ contains a subsequence of type $\alpha^i_j$ then that subsequence is entirely contained in either $|\sigma_1|+2\ \sigma_1^{+1} 1$ or $\sigma_2^{+|\sigma_1|+2}$.
As a result, the set of involutions which avoid 3412 and $T$ and which do not begin with 1 can be partitioned into sets $B_{i_1,\ldots,i_m}$, where $B_{i_1,\ldots,i_m}$ is the set of such involutions in which $\sigma_1$ contains $T^{left}_{i_1-1,\ldots,i_m-1}$ but avoids $T^{left}_{i_1,\ldots,i_m}$ and $\sigma_2$ avoids $T^{right}_{i_1,\ldots,i_m}$.
In view of Lemma \ref{lem:Y}, the generating function for $B_{i_1,\ldots,i_m}$ is $\left(\sum\limits_{Y \subseteq \{1,2,\ldots,m\}} (-1)^{|Y|} F_{T_Y}(x)\right) F_{T^{right}_{i_1,\ldots,i_m}}(x)$.
It follows that the generating function for $A_3$ is 
$$x^2 \sum_{i_1,\ldots,i_m = 1}^{k_1,\ldots,k_m} \left( \sum_{Y \subseteq \{1,2,\ldots,m\}} (-1)^{|Y|} F_{T_Y}(x) \right) F_{T^{right}_{i_1,\ldots,i_m}}(x).$$

Add the generating functions for $A_1$, $A_2$ and $A_3$ to obtain (\ref{eqn:FTrecurrence}).
\end{proof}

The case of (\ref{eqn:FTrecurrence}) in which $|T| = 1$ will prove useful, so we single it out here.

\begin{corollary}
Suppose $\pi = \alpha_1|\ldots|\alpha_k$ is a permutation, where $\alpha_1,\ldots,\alpha_k$ are complete.
Then
\begin{equation}
\label{eqn:Fpirecurrence}
F_\pi(x) = 1 + x F_\beta(x) + x^2 \sum_{i=1}^k \left( F_{\overline{\alpha_1|\ldots |\alpha_i}}(x) - F_{\overline{\alpha_1|\ldots|\alpha_{i-1}}}(x)\right) F_{\alpha_i|\ldots|\alpha_k}(x).
\end{equation}
Here $\beta = \pi$ if $\alpha_1 \neq 1$ and $\beta = \alpha_2|\ldots|\alpha_k$ if $\alpha_1 = 1$.
\end{corollary}
\begin{proof}
Set $T = \{\pi\}$ in Theorem \ref{thm:FTrecurrence}.
\end{proof}

Observe that if $\overline{\pi} \neq \pi$ or $\pi$ is not complete then (\ref{eqn:Fpirecurrence}) allows one to express $F_\pi(x)$ in terms of $F_\sigma(x)$ for various $\sigma$ with $|\sigma| < |\pi|$.
If $\overline{\pi} = \pi$ and $\pi$ is complete then one can solve (\ref{eqn:Fpirecurrence}) to find that $F_\pi(x) = \sum\limits_{n=0}^\infty |I_n(3412)| x^n$.
Since $I_n(3412, \pi) \subseteq I_n(3412)$, we have the following result, which can also be shown with a routine induction argument.

\begin{proposition}
\label{prop:completebar}
If $\overline{\pi} = \pi$ and $\pi$ is complete then $I_n(3412, \pi) = I_n(3412)$ for all $n \ge 0$.
\end{proposition}

\section{Generating Functions Involving Chebyshev Polynomials}

In this section we use (\ref{eqn:Fpirecurrence}) to find $F_\pi(x)$ for various $\pi$.
In each case we express $F_\pi(x)$ in terms of Chebyshev polynomials of the second kind.
We begin with $F_{k\ldots 21}(x)$.

\begin{theorem}
\label{thm:F321}
For all $k \ge 1$, the following hold.
\begin{equation}
\label{eqn:F2k}
F_{2k \ldots 21}(x) = \frac{U_{k-1}\left(\frac{1-x}{2x}\right)}{x U_k\left(\frac{1-x}{2x}\right)}
\end{equation}
\begin{equation}
\label{eqn:F2k-1}
F_{2k-1 \ldots 21}(x) = \frac{U_{k-1}\left(\frac{1-x}{2x}\right) + U_{k-2}\left(\frac{1-x}{2x}\right)}{x \left( U_{k}\left(\frac{1-x}{2x}\right) + U_{k-1}\left(\frac{1-x}{2x}\right)\right)}
\end{equation}
\end{theorem}
\begin{proof}
To prove (\ref{eqn:F2k}), we argue by induction on $k$.
It is routine to verify that the result holds for $k=1$, so assume it holds for $k-1$.
Set $\pi = 2k \ldots 21$ in (\ref{eqn:Fpirecurrence}) and solve the resulting equation for $F_{2k \ldots 21}(x)$ to obtain
$$F_{2k \ldots 21}(x) = \frac{1}{1 - x - x^2 F_{2k-2 \ldots 21}(x)}.$$
Now use induction to eliminate $F_{2k-2 \ldots 21}(x)$ and simplify the result to obtain
$$F_{2k \ldots 21}(x) = \frac{U_{k-1}\left(\frac{1-x}{2x}\right)}{(1-x) U_{k-1}\left( \frac{1-x}{2x}\right) - x U_{k-2}\left( \frac{1-x}{2x} \right)}.$$
Finally, use (\ref{eqn:ourChebyrecurrence}) to simplify the denominator and obtain (\ref{eqn:F2k}).

The proof of (\ref{eqn:F2k-1}) is similar to the proof of (\ref{eqn:F2k}).
\end{proof}

\noindent
{\bf Remark}
Lines (\ref{eqn:F2k}) and (\ref{eqn:F2k-1}) can also be obtained using $\varphi$ and the methods of Section \ref{sec:rk21}.

\medskip

Arguing by induction as in the proof of (\ref{eqn:F2k}) and using Simion and Schmidt's result \cite[Prop. 6]{SimionSchmidt} that $F_{312}(x) = F_{231}(x) = \frac{1-x}{1-2x}$, we obtain the following result.

\begin{theorem}
\label{thm:F231312}
For all $k \ge 3$, the following holds.
\begin{equation}
\label{eqn:Fk312}
F_{k+1\ k \ldots 4231}(x) = F_{k\ldots 312}(x) = \frac{U_{k-2}\left(\frac{1-x}{2x}\right)}{x U_{k-1}\left(\frac{1-x}{2x}\right)}
\end{equation}
\end{theorem}

The fact that $F_{k+1\ k\ldots 4231}(x) = F_{k\ldots 312}(x)$ implies that $I_n(3412, k+1\ k\ldots 4231)$ and $I_n(3412, k\ldots 312)$ have the same cardinality.
But more is true:  as we show next, $I_n(3412, k+1\ k\ldots 4231)$ and $I_n(3412, k\ldots 312)$ are the same set.

\begin{theorem}
\label{thm:kk1equalassets}
Suppose $\sigma$ is a nonempty permutation which does not end with 1.
Then for all $n \ge 0$,
\begin{equation}
\label{eqn:kk1equalassets}
I_n(3412, |\sigma|+1\ \sigma) = I_n(3412, |\sigma|+2\ \sigma^{+1}\ 1).
\end{equation}
\end{theorem}
\begin{proof}
It is clear that if $\pi$ avoids $|\sigma|+1\ \sigma$ then $\pi$ avoids $|\sigma|+2\ \sigma^{+1}\ 1$, so $I_n(3412, |\sigma|+1\ \sigma) \subseteq I_n(3412, |\sigma|+2\ \sigma^{+1}\ 1)$.

Now suppose $\pi$ contains $|\sigma|+1\ \sigma$;  we show $\pi$ contains $|\sigma|+2\ \sigma^{+1}\ 1$.
We argue by induction on $|\pi|$.

The result is vacuously true if $|\pi| \le |\sigma|$, and it follows for $|\pi| = |\sigma|+1$ since $|\sigma|+1\ \sigma$ does not end with 1, and so is not an involution.
Now suppose the result holds for all involutions in $I(3412)$ of length less than $|\pi|$.
If $\pi$ begins with 1 then the result follows by induction, in view of Proposition \ref{prop:permmap}(i).
If $\pi$ begins with $j \ge 2$ then by Proposition \ref{prop:permmap}(ii) there exist $\pi_1 \in I_{j-2}(3412)$ and $\pi_2 \in I_{n-j}(3412)$ such that $\pi = j\ \pi_1^{+1}\ 1\ \pi_2^{+j}$.
Since $|\sigma|+1\ \sigma$ begins with its greatest entry, it is entirely contained in either $j\ \pi_1^{+1}\ 1$ or $\pi_2^{+j}$.
In the latter case the result follows by induction.
In the former, $|\sigma|+1\ \sigma$ is contained in $j \pi_1$, since $\sigma$ does not end with 1.
It follows that $|\sigma|+2\ \sigma^{+1}\ 1$ is contained in $j\ \pi_1\ 1$, as desired.
\end{proof}

\begin{corollary}
For all $n \ge 0$ and all $k \ge 4$, the sets $I_n(3412, k \ldots 312)$ and $I_n(3412, k+1 \ldots 4231)$ are equal.
\end{corollary}
\begin{proof}
Set $\sigma = k-1 \ldots 312$ in (\ref{eqn:kk1equalassets}).
\end{proof}

We now turn our attention to $F_{k\ldots 4213}(x)$, $F_{k\ldots 4132}(x)$, and $F_{k\ldots 4123}(x)$.

\begin{theorem}
\label{thm:F213132}
For all $k \ge 3$, the following holds.
\begin{equation}
\label{eqn:F213132}
F_{k\ldots 4213}(x) = F_{k\ldots 4132}(x) = \frac{U_{k-2}\left( \frac{1-x}{2x}\right) + U_{k-3}\left(\frac{1-x}{2x}\right)}{x \left( U_{k-1}\left(\frac{1-x}{2x}\right) + U_{k-2}\left(\frac{1-x}{2x}\right)\right)}
\end{equation}
\end{theorem}
\begin{proof}
First observe that if we apply the reverse complement map to $k\ldots 4213$ and take the inverse of the result we obtain $k\ldots 4132$, so $F_{k\ldots 4213}(x) = F_{k\ldots 4132}(x)$.

To show that $F_{k\ldots 4132}(x)$ is equal to the quantity on the right, first observe that the result holds for $k=3$ by \cite[Ex. 2.18]{MansourGuibert}.
Now argue by induction on $k$, using (\ref{eqn:Fpirecurrence}) and (\ref{eqn:ourChebyrecurrence}).
\end{proof}

\begin{theorem}
\label{thm:F123}
For all $k \ge 3$, the following holds.
$$F_{k\ldots 4123}(x) = \frac{(1-x+x^3) U_{k-3}\left(\frac{1-x}{2x}\right) + (x-1) x U_{k-4}\left(\frac{1-x}{2x}\right)}{(1-x+x^3) x U_{k-2}\left(\frac{1-x}{2x}\right) + (x-1)x^2U_{k-3}\left(\frac{1-x}{2x}\right)}$$
\end{theorem}
\begin{proof}
This is similar to the second half of the proof of Theorem \ref{thm:F213132}.
\end{proof}

Recall that 213 and 123 are examples of layered permutations, which are defined as follows.

\begin{definition}
Fix $n \ge 1$ and let $l_1, l_2, \ldots, l_m$ denote a sequence such that $l_i \ge 1$ for $1 \le i \le m$ and $\sum\limits_{i=1}^m l_i = n$.
We write $[l_1, l_2, \ldots, l_m]$ to denote the permutation given by
$$[l_1, l_2, \ldots, l_m] = l_1, l_1-1, \ldots, 1, l_2 + l_1, l_2+l_1-1, \ldots, l_1+1, \ldots, n, n-1, \ldots, n-l_m+1.$$
We call a permutation {\em layered} whenever it has the form $[l_1,\ldots,l_m]$ for some sequence $l_1, \ldots, l_m$.
\end{definition}

Observe that if $m \ge 2$ then $\overline{[l_1, \ldots, l_m]} = [l_1,\ldots,l_m]$.
In view of (\ref{eqn:Fpirecurrence}), (\ref{eqn:F2k}), and (\ref{eqn:F2k-1}), the generating function $F_{[l_1, \ldots, l_m]}(x)$ can be expressed in terms of Chebyshev polynomials of the second kind for any layered permutation $[l_1,\ldots, l_m]$.
For example, when $m = 2$ we have the following result.

\begin{theorem}
\label{thm:2layered}
For all $k, l \ge 1$ we have
\begin{equation}
\label{eqn:2layered}
F_{[k,l]}(x) = F_{k+l\ k+l-1\ldots 21}(x).
\end{equation}
\end{theorem}

To prove Theorem \ref{thm:2layered}, we need the following well-known result concerning Chebyshev polynomials.

\begin{lemma}
\label{lem:Chebysum}
For all $k, l \ge -1$ and all $w \ge 0$ we have
\begin{equation}
\label{eqn:Chebysum}
U_{k+w} U_{l+w} - U_k U_l = U_{w-1} U_{k+l+w+1}.
\end{equation}
Here we abbreviate $U_* = U_*\left(\frac{1-x}{2x}\right)$.
\end{lemma}
\begin{proof}
We argue by induction on $k+w$.

Line (\ref{eqn:Chebysum}) is immediate for $k+w = -1$ and $k+w = 0$ so suppose $k+w > 0$ and (\ref{eqn:Chebysum}) holds for $k+w-1$ and $k+w-2$.
Divide (\ref{eqn:ourChebyrecurrence}) by $x$ and use the result to eliminate $U_{k+w}$ and $U_k$ on the left side of (\ref{eqn:Chebysum}).
Now use induction and (\ref{eqn:ourChebyrecurrence}) twice each to obtain (\ref{eqn:Chebysum}).
\end{proof}

We now prove Theorem \ref{thm:2layered}.

\noindent
{\em Proof of Theorem \ref{thm:2layered}}.
We consider six cases:  $k = 1$ and $l$ is even, $k=1$ and $l$ is odd, $k$ and $l$ are both even, $k$ is even and $l$ is odd, $k >1$ is odd and $l$ is even, and $k>1$ and $l$ are both odd.
All six cases are similar, so we only give the details for the case in which $k>1$ and $l$ are both odd.

Set $\pi = [2k+1, 2l-1]$ in (\ref{eqn:Fpirecurrence}) and solve the resulting equation for $F_{[2k+1,2l-1]}(x)$ to obtain
$$F_{[2k+1,2l-1]}(x) = \frac{1-x^2 F_{[2k-1]}(x) F_{[2l-1]}(x)}{1-x-x^2 F_{[2k-1]}(x) - x^2 F_{[2l-1]}(x)}.$$
Now use (\ref{eqn:F2k-1}) to eliminate $F_{[2k-1]}(x)$ and $F_{[2l-1]}(x)$, clear denominators, and use (\ref{eqn:ourChebyrecurrence}) to simplify the denominator of the result, obtaining
$$F_{[2k+1,2l-1]}(x) = \frac{(U_k+U_{k-1})(U_l+U_{l-1}) - (U_{k-1} + U_{k-2})(U_{l-1}+U_{l-2})}{x(U_{k+1}+ U_k)(U_l+U_{l-1}) - x(U_k+U_{k-1})(U_l + U_{l-1})}.$$
Here we abbreviate $U_* = U_*\left(\frac{1-x}{2x}\right)$.
Next use (\ref{eqn:Chebysum}) with $w=1$ four times each in the numerator and denominator to obtain
$$F_{[2k+1,2l-1]}(x) = \frac{U_{k+l} + 2 U_{k+l-1} + U_{k+l-2}}{x\left( U_{k+l+1} + 2 U_{k+l} + U_{k+l-1}\right)}.$$
Rearranging (\ref{eqn:ourChebyrecurrence}) we find that $U_n + U_{n-1} + U_{n-2} = \frac{1}{x} U_{n-1}$ for all $n \ge 1$.
Use this to simplify our last expression for $F_{[2k+1,2l-1]}(x)$, obtaining the right side of (\ref{eqn:F2k}) with $k$ replaced by $k+l$, as desired.
$\Box$

When $m \ge 3$ the generating function $F_{[l_1,\ldots,l_m]}(x)$ does not reduce quite as nicely as it does when $m=2$.
For example, using the same techniques as in the proof of Theorem \ref{thm:2layered} one can prove that for all $k_1, k_2, k_3 \ge 1$,
\begin{equation}
\label{eqn:k1k2k3}
F_{[2k_1,2k_2,2k_3]}(x) = \frac{U_{k_1+k_2+k_3} U_{k_1+k_2+k_3-1} + U_{k_1+k_2-1} U_{k_1+k_3-1} U_{k_2+k_3-1}}{x U_{k_1+k_2} U_{k_1+k_3} U_{k_2+k_3}}.
\end{equation}
Nevertheless, (\ref{eqn:2layered}) and (\ref{eqn:k1k2k3}) suggest the following conjecture.

\begin{conjecture}
\label{conj:symmetry}
For all $m \ge 1$ and all $l_1,\ldots, l_m  \ge 1$, the generating function $F_{[l_1,\ldots,l_m]}(x)$ is symmetric in $l_1,\ldots,l_m$.
\end{conjecture}

This conjecture has been verified for $m \le 4$ and $l_i \le 20$, as well as for $m = 5$ and $l_i \le 11$, using a Maple program.

\section{Involutions Which Avoid 3412 and Several Other Patterns}

Let $T$ denote a set of permutations, each of which begins with its largest element and ends with 1.
In this section we use (\ref{eqn:FTrecurrence}) to find the generating function for those involutions which avoid 3412 and every element of $T$ in terms of the generating function for those involutions which avoid 3412 and the permutations obtained by removing the first and last entry from each element of $T$.
We then use this recurrence to enumerate involutions which avoid 3412 and various sets of permutations.
We begin by setting some notation.

\begin{definition}
For any set $T$ of permutations, we write $k T 1$ to denote the set of permutations obtained by replacing each $\pi \in T$ with $|\pi|+2\ \pi^{+1}\ 1$.
For all $i \ge 2$, we also write $k^i T 1^i$ to denote the set obtained by applying this operation to $T$ $i$ times.
\end{definition}

Our main result, which we prove next, gives $F_{kT1}(x)$ in terms of $F_T(x)$.

\begin{proposition}
\label{prop:kT1}
For any set $T$ of permutations, we have
\begin{equation}
\label{eqn:kT1}
F_{kT1}(x) = \frac{1}{1-x-x^2 F_T(x)}.
\end{equation}
\end{proposition}
\begin{proof}
Replace $T$ with $kT1$ in (\ref{eqn:FTrecurrence}), use the fact that every permutation in $kT1$ is complete and $\beta(kT1) = kT1$ to simplify the result, and solve for $F_{kT1}(x)$.
\end{proof}

Proposition \ref{prop:kT1} leads to an enumeration involving a generalization of the Fibonacci numbers.

\begin{corollary}
For all $n \ge 0$ and all $k \ge 1$ we have
$$|I_n(3412, 4231, k+2\ k+1\ \ldots 4321)| = F_{k+1,n+1},$$
where $F_{k,n}$ is the $k$-generalized Fibonacci number, defined by $F_{k,n} = 0$ for $n \le 0$, $F_{k,1} = 1$, and $F_{k,n} = \sum\limits_{i=1}^k F_{k,n-i}$ for all $n \ge 1$.
\end{corollary}
\begin{proof}
Set $T = \{k\ldots 21, 12\}$ in (\ref{eqn:kT1}) and observe that $F_T(x) = \sum\limits_{i=0}^{k-1} x^i$.
\end{proof}

If $F_T(x)$ is a rational function then $F_{k^i T 1^i}(x)$ can be expressed in terms of Chebyshev polynomials of the second kind, as we show next.

\begin{corollary}
Suppose $F_T(x) = \frac{f_0(x)}{f_1(x)}$ for polynomials $f_0$ and $f_1$.
Then for all $i \ge 1$,
\begin{equation}
\label{eqn:FkiT1i}
F_{k^i T 1^i}(x) = \frac{f_1(x) U_{i-1}\left(\frac{1-x}{2x}\right) - x f_0(x) U_{i-2}\left(\frac{1-x}{2x}\right)}{x f_1(x) U_i\left(\frac{1-x}{2x}\right) - x^2 f_0(x) U_{i-1}\left(\frac{1-x}{2x}\right)}.
\end{equation}
\end{corollary}
\begin{proof}
We argue by induction on $i$.

When $i = 1$ the result is immediate from (\ref{eqn:kT1}), so we assume the result holds for $i-1$.
Replace $T$ with $k^{i-1} T 1^{i-1}$ in (\ref{eqn:kT1}), use induction to eliminate $F_{k^{i-1} T 1^{i-1}}(x)$ on the right, and use (\ref{eqn:ourChebyrecurrence}) to simplify the result and obtain (\ref{eqn:FkiT1i}).
\end{proof}

\noindent
{\bf Remark}
For any set $T$ of permutations, let $kT$ denote the set of permutations obtained by replacing each permutation $\pi \in T$ with $|\pi|+1\ \pi$.
If no permutation in $T$ ends with 1, then by Theorem \ref{thm:kk1equalassets} all of the results in this section hold when $kT1$ is replaced with $kT$.

\section{Enumerations}

In this section we give enumerations of $I_n(3412, \sigma)$ for various $\sigma$;  our results are listed in the tables below.
We include references when they exist;  if no reference is given then the enumeration follows from (\ref{eqn:Fpirecurrence}).
We place $\sigma_1$ and $\sigma_2$ on the same line whenever $\sigma_2$ can be obtained from $\sigma_1$ using the inverse and reverse complement maps, since in this case $|I_n(3412,\sigma_1)| = |I_n(3412,\sigma_2)|$ for all $n\ge 0$.

Several of our enumerations are expressed in terms of Motzkin, Fibonacci, or Pell numbers.
We write $M_n$ to denote the $n$th Motzkin number, which may be defined by $M_0 = 1$ and $M_n = M_{n-1} + \sum\limits_{i=2}^n M_{i-2} M_{n-i}$ for all $n \ge 1$.
We write $F_n$ to denote the $n$th Fibonacci number, which may be defined by $F_0 = 0$, $F_1 = 1$, and $F_n = F_{n-1} + F_{n-2}$ for all $n \ge 2$.
We write $P_n$ to denote the $n$th Pell number, which may be defined by $P_0 = 0$, $P_1 = 1$, and $P_n = 2 P_{n-1} + P_{n-2}$ for all $n \ge 2$.

\bigskip

\renewcommand{\baselinestretch}{1.3}
\small
\normalsize

\begin{center}
\begin{tabular}{|c|c|c|c|c|c|}
\hline
\multicolumn{4}{|c|}{$\sigma$} & $|I_n(3412,\sigma)|$ & Reference \\
\hline
\hline
\multicolumn{4}{|c|}{123} & $\frac{2n^2 + 7 + (-1)^n}{8}$ & Theorem \ref{thm:F123} \\
\hline
\multicolumn{2}{|c|}{132} & \multicolumn{2}{|c|}{213} & & \cite[Ex. 2.18]{MansourGuibert} \\
\cline{1-4}
\cline{6-6}
\multicolumn{4}{|c|}{321} & \raisebox{1.5ex}[0pt]{$F_{n+1}$} & Theorem \ref{thm:F321} \\
\hline
\multicolumn{2}{|c|}{231} & \multicolumn{2}{|c|}{312} & & \cite[Prop. 6]{SimionSchmidt} \\
\cline{1-4}
\cline{6-6}
\multicolumn{2}{|c|}{1432} & \multicolumn{2}{|c|}{3214} & & Theorem \ref{thm:2layered} \\
\cline{1-4}
\cline{6-6}
\multicolumn{4}{|c|}{2143} & $2^{n-1}$ & Theorem \ref{thm:2layered} \\
\cline{1-4}
\cline{6-6}
\multicolumn{4}{|c|}{4231} & & Theorem \ref{thm:F231312} \\
\cline{1-4}
\cline{6-6}
\multicolumn{4}{|c|}{4321} & & Theorem \ref{thm:F321} \\
\hline
\end{tabular}

\bigskip
\bigskip

\begin{tabular}{|c|c|c|c|c|c|}
\hline
\multicolumn{4}{|c|}{$\sigma$} & $|I_n(3412,\sigma)|$ & Reference \\
\hline
\hline
\multicolumn{4}{|c|}{1234} & $\frac{2n^4-4n^3+28n^2-2n+81-6n(-1)^n+15(-1)^n}{96}$ &  \\
\hline
\multicolumn{2}{|c|}{1243} & \multicolumn{2}{|c|}{2134} & & \\
\cline{1-4}
\cline{6-6}
\multicolumn{4}{|c|}{1324} & \raisebox{1.5ex}[0pt]{$\frac{n}{5} F_{n+2} + \frac{n}{5} F_n-\frac{3}{5} F_n + 1$} &  \\
\hline
1342 & 1423 & 2314 & 3124 & & \\
\cline{1-4}
\cline{6-6}
2431 & 4132 & 3241 & 4213 & & \\
\cline{1-4}
\cline{6-6}
\multicolumn{2}{|c|}{21543} & \multicolumn{2}{|c|}{32154} & & Theorem \ref{thm:2layered} \\
\cline{1-4}
\cline{6-6}
\multicolumn{2}{|c|}{43215} & \multicolumn{2}{|c|}{15432} & $g.f. = \frac{1-x-x^2}{1-2x-x^2+x^3}$ & Theorem \ref{thm:2layered} \\
\cline{1-4}
\cline{6-6}
\multicolumn{2}{|c|}{53241} & \multicolumn{2}{|c|}{52431} & & Proposition \ref{prop:kT1} \\
\cline{1-4}
\cline{6-6}
\multicolumn{2}{|c|}{42315} & \multicolumn{2}{|c|}{15342} & & \\
\cline{1-4}
\cline{6-6}
\multicolumn{4}{|c|}{54321} & & Theorem \ref{thm:F321} \\
\hline
\multicolumn{4}{|c|}{52341} & & Proposition \ref{prop:kT1} \\
\cline{1-4}
\cline{6-6}
\multicolumn{2}{|c|}{2341} & \multicolumn{2}{|c|}{4123} & \raisebox{1.5ex}[0pt]{$g.f. = \frac{(1-x)^2 (1-x^2)}{1-3x+x^2+3x^3-3x^4}$} & Theorem \ref{thm:F123} \\
\hline
\multicolumn{2}{|c|}{2413} & \multicolumn{2}{|c|}{3142} & & Proposition \ref{prop:completebar} \\
\cline{1-4}
\cline{6-6}
\multicolumn{4}{|c|}{3412} & \raisebox{1.5ex}[0pt]{$M_n$} & \cite[Rem. 4.28]{G} \\
\hline
\multicolumn{2}{|c|}{3421} & \multicolumn{2}{|c|}{4312} & & Theorem \ref{thm:F231312} \\
\cline{1-4}
\cline{6-6}
\multicolumn{2}{|c|}{32541} & \multicolumn{2}{|c|}{52143} & & \\
\cline{1-4}
\cline{6-6}
51432 & 43251 & 25431 & 53214 & & \\
\cline{1-4}
\cline{6-6}
14352 & 15324 & 41325 & 24315 & $\frac{P_n+P_{n-1}+1}{2}$ & \\
\cline{1-4}
\cline{6-6}
21534 & 23154 & 21453 & 31254 & & \\
\cline{1-4}
\cline{6-6}
13542 & 15243 & 42135 & 32415 & & \\
\cline{1-4}
\cline{6-6}
\multicolumn{2}{|c|}{54231} & \multicolumn{2}{|c|}{53421} & & Proposition \ref{prop:kT1} \\
\hline
\end{tabular}

\bigskip
\bigskip

\begin{tabular}{|c|c|c|c|c|c|}
\hline
\multicolumn{4}{|c|}{$\sigma$} & $|I_n(3412,\sigma)|$ & Reference \\
\hline
\hline
\multicolumn{2}{|c|}{51324} & \multicolumn{2}{|c|}{24351} & & \\
\cline{1-4}
\cline{6-6}
32451 & 51243 & 52134 & 23541 & \raisebox{1.5ex}[0pt]{$g.f. = \frac{(1-x)(1-x-x^2)^2}{1-4x+3x^2+4x^3-4x^4-2x^5+x^6}$} & \\
\hline
\multicolumn{4}{|c|}{21354} & & \\
\cline{1-4}
\cline{6-6}
\multicolumn{2}{|c|}{13254} & \multicolumn{2}{|c|}{21435} & \raisebox{1.5ex}[0pt]{$3 \cdot 2^{n-1} - \frac{2 F_{n+3} + n F_{n+2} + F_{n+1}+nF_n}{5}$} & \\
\hline
24531 & 51423 & 53124 & 34251 & & \\
\cline{1-4}
\cline{6-6}
43521 & 54213 & 54132 & 35421 & & \\
\cline{1-4}
\cline{6-6}
42351 & 52314 & 51342 & 25341 & \raisebox{1.5ex}[0pt]{$g.f. = \frac{1-2x-x^2+x^3}{1-3x+3x^3}$} & \\
\cline{1-4}
\cline{6-6}
34215 & 15423 & 43125 & 14532 & & \\
\hline
\multicolumn{2}{|c|}{12543} & \multicolumn{2}{|c|}{32145} & & \\
\cline{1-4}
\cline{6-6}
\multicolumn{4}{|c|}{14325} & \raisebox{1.5ex}[0pt]{$\frac{(3n+7) 2^{n-2}}{9} + \frac12 + \frac{1}{18} (-1)^n$} & \\
\hline
15234 & 13452 & 23415 & 41235 & $g.f. = \frac{1-2x-x^2+3x^3-x^4-x^5}{1-3x+5x^3-3x^4-2x^5+x^6}$ & \\
\hline
12534 & 12453 & 23145 & 31245 & & \\
\cline{1-4}
\cline{6-6}
\multicolumn{2}{|c|}{13425} & \multicolumn{2}{|c|}{14235} & \raisebox{1.5ex}[0pt]{$g.f. = \frac{1-4x+3x^2+4x^3-4x^4+x^6}{(1-x)(1-2x-x^2+x^3)^2}$} & \\
\hline
\multicolumn{2}{|c|}{21345} & \multicolumn{2}{|c|}{12354} & & \\
\cline{1-4}
\cline{6-6}
\multicolumn{2}{|c|}{12435} & \multicolumn{2}{|c|}{13245} & \raisebox{1.5ex}[0pt]{$\frac{(5n^2-3n-100)F_{n+1}}{50} +\frac{(12n-38) F_n}{25} + \frac{2n^2+8n+23+(-1)^n}{8}$} & \\
\hline
\multicolumn{2}{|c|}{54312} & \multicolumn{2}{|c|}{45321} & $\frac12 (F_{n+1} + F_{2n} - F_{2n-2})$ & Theorem \ref{thm:F231312} \\
\hline
\multicolumn{2}{|c|}{34521} & \multicolumn{2}{|c|}{54123} & $g.f. = \frac{1-3x+x^2+3x^3-3x^4}{(1-x)(1-3x+x^2+3x^3-3x^4)}$ & Theorem \ref{thm:F123} \\
\hline
\multicolumn{4}{|c|}{12345} & $\frac{2n^6-12n^5+86n^4-168n^3+731n^2-54n+1917 +(-1)^n (45n^2-234n+387)}{2304}$ & \\
\hline
\multicolumn{2}{|c|}{51234} & \multicolumn{2}{|c|}{23451} & $g.f. = \frac{(1-x)^5 (1+x)^2}{1-4x+3x^2+6x^3-10x^4+6x^6-4x^7}$ & \\
\hline
\end{tabular}
\end{center}

\renewcommand{\baselinestretch}{1}
\small
\normalsize

\section{Directions for Future Research}

\begin{enumerate}
\item
For any permutation $\pi \neq k\ldots 21$, let $G_\pi(x)$ denote the generating function for the set of involutions which avoid 3412 and contain exactly one subsequence of type $\pi$.
Compute $G_\pi(x)$ or find a recursive formula for $G_\pi(x)$.
More generally, find the generating function for the set of involutions which avoid 3412 and contain exactly $r$ subsequences of type $\pi$.
\item
In view of (\ref{eqn:F2k}) and (\ref{eqn:Fk312}) we have
$$|I_n(3412,2k\ldots 21)| = |I_n(3412, k+2\ k+1\ldots 4231)|$$
for all $n\ge 0$.
Similarly, in view of (\ref{eqn:F2k-1}) and (\ref{eqn:F213132}) we have
$$|I_n(3412,2k+1\ldots 21)| = |I_n(3412,k+2\ k+1\ldots 4132)|$$
for all $n\ge 0$.
Give combinatorial proofs of these identities.
\item
In view of (\ref{eqn:2layered}), we have
$$|I_n(3412, [k,l])| = |I_n(3412, k+l\ k+l-1\ldots 21)|$$
for all $k,l \ge 1$ and all $n \ge 0$.
Give a combinatorial proof of this identity.
\item
Prove Conjecture \ref{conj:symmetry}, which says that $F_{[l_1,\ldots,l_m]}(x)$ is symmetric in $l_1,\ldots, l_m$.
\end{enumerate}

\begin{Large}
\noindent
{\bf Acknowledgement}
\end{Large}

\medskip

The author thanks Toufik Mansour for several helpful comments and suggestions on this paper.


\end{document}